%%%%%%%%%%%%%%  Geometry and Topology Monogrpahs: m4-20.tex  %%%%%%%%
%%%%        
%%%%    On the potential functions for the hyperbolic structures ...
%%%%             
%%%%                         Yoshiyuki Yokota
%%%%                          
%%%%                  Published in Volume 4(2002) 303-311
%%%%
%%%%                    Publication date 13 October 2002
%%%%
%%%%                      This is a plain TeX file
%%%%
%%%%
%%%%%%%%%%%%%%%%%%                                   %%%%%%%%%%%%%%%%%%%
%%%%%%%%%%%%%%%%%%%%%%%%%%%%%%%%%%%%%%%%%%%%%%%%%%%%%%%%%%%%%%%
%%%%%%%%%%%             gtmacros.tex            %%%%%%%%%%%%%%%
%%%%%%%%%%%             version 1.6             %%%%%%%%%%%%%%% 
%
%                       Colin Rourke   
%
%
%    These macros are recommended for use by authors submitting articles   
%    to Geometry and Topology or to Algebraic and Geometric Topology.  
%    They are intended to be used with plain TeX. Each macro is described 
%    briefly to make it clear how to use it (or to modify it to achieve
%    different results).  If you modify this file then please change its
%    name.  If you modify this file and use the modified file to 
%    format an article for submission to Geometry and Topology or
%    Algebraic and Geometric Topology, then please paste the modified
%    file into your main TeX file.  Do not submit it as a separate file.
%      
%    Instructions on using these macros are also given in  gtmacins.tex  
%    or  gtmacins.ps  or .pdf  available on the gt www pages or by 
%    anonymous ftp from the gt/info/macros directory.
%
%
\magnification=\magstephalf      % Sets default point size to 11pt.
%
%  Basic layout parameters :
%
\vsize=7.5truein                 % Sets text height to 7.5 inches.
\hsize=5.2truein                 % Sets text width to 5.2 inches.
\newskip\stdskip                 % standard vertical space
\stdskip=6pt plus3pt minus3pt    % (slightly more stretchy
\medskipamount=\stdskip          % than the usual \medskip)
\parindent=0pt                   % Paragraphs are non-indented with
\parskip=\stdskip                % a little space between paragraphs. 
\abovedisplayskip=\stdskip       %  Reduces the space
\belowdisplayskip=\stdskip       %  around displays.
\mathsurround=0.75pt             % Gives a little extra space around maths.
\overfullrule=0pt                %  Prevents black boxes
%
%   The following macro is for principal paragraph breaks ie
%   a paragraph break with a slightly larger space :
%
\def\ppar{\par\goodbreak\vskip 8pt plus 4pt minus 4pt}     
%
%  The standard horizontal space for theorems, labels etc :
%
\def\stdspace{\hskip 0.75em plus 0.15em\ignorespaces}
\let\qua\stdspace % useful abbreviation (3/4 of a quad)
%
%%%%%%%%%%%%%%            FONT MACROS            %%%%%%%%%%%%
%
%           The following font macros define the AMS symbol 
%           and Euler-Fraktal fonts for use in text and
%           mathematics with appropriate size changes.
%           They also define two new control sequences  
%           \small  and  \large  (similar to those built
%           into LaTeX) which change the size of all fonts 
%           both in text and maths.  \small  is 10% smaller 
%           than normal and  \large  30% bigger.  The strange
%           size of the \large text fonts (10pt scaled 1315)
%           is because these macros are intended to be used
%           at \magstephalf.  The result is 10pt scaled 1440
%           (\magstep2) which is a standard font size.  If
%           you are borrowing these macros to use them at
%           another basic  \magnification, then you will
%           probably need to change 1315 to 1200 in the eleven
%           places marked ** below.  \large  will then be
%           20% bigger than normal.  Note that at \magstephalf
%           all the fonts come out roughly one point larger
%           than their size as defined in these macros.
%
%           The size-changing macros are based on Knuth's
%           \ninepoint and \eightpoint macros.
%
%
%    The macros are laid out in a way which makes it clear how to
%    add futher fonts (or delete unavailable ones) and how to add
%    further size changes.
%
%    First comes a definition of  \hexnumber  which is needed to
%    refer to font families whose family number is not known :
%
\def\hexnumber#1{\ifcase#1 0\or 1\or 2\or 3\or 4\or 5\or 6\or 7\or 8\or
 9\or A\or B\or C\or D\or E\or F\fi}
%
%     Next we define the AMS symbol-a fonts at 13,10,9,7,6,5 points
%
\font\thirtnmsa=msam10 scaled 1315    %%% **  see note above 
\font\tenmsa=msam10          \font\ninemsa=msam9
\font\sevenmsa=msam7         \font\sixmsa=msam6
\font\fivemsa=msam5
%%%%%%  (add further sizes here if you need them)
%
%    and the standard size family for these fonts
%
\newfam\msafam                  \textfont\msafam=\tenmsa
\scriptfont\msafam=\sevenmsa    \scriptscriptfont\msafam=\fivemsa
\edef\hexa{\hexnumber\msafam}        %  The msa family is  \fam\hexa
\def\msa{\fam\msafam\tenmsa}         %  \msa  switches to this family
%
%    Repeat these steps for the AMS symbol-b fonts
%
\font\thirtnmsb=msbm10 scaled 1315   %%%  ** see note above
\font\tenmsb=msbm10      \font\ninemsb=msbm9
\font\sevenmsb=msbm7     \font\sixmsb=msbm6
\font\fivemsb=msbm5
%%%%%%  (add further sizes here if you need them)
%
\newfam\msbfam                   \textfont\msbfam=\tenmsb       
\scriptfont\msbfam=\sevenmsb     \scriptscriptfont\msbfam=\fivemsb
\edef\hexb{\hexnumber\msbfam}    %  The msb family is \fam\hexb  
\def\msb{\fam\msbfam\tenmsb}     %  \msb switches to this family
%
%        Repeat for the Euler-Fraktal fonts 
%
\font\thirtneufm=eufm10 scaled 1315   %%% **  see note above 
\font\teneufm=eufm10                 \font\nineeufm=eufm9
\font\seveneufm=eufm7                \font\sixeufm=eufm6
\font\fiveeufm=eufm5
%%%%%%  (add further sizes here if you need them)
%
\newfam\eufmfam                    \textfont\eufmfam=\teneufm
\scriptfont\eufmfam=\seveneufm     \scriptscriptfont\eufmfam=\fiveeufm
\edef\hexf{\hexnumber\eufmfam}      % The Euler-Fraktal family is
\def\frak{\fam\eufmfam\teneufm}     % \fam\hexf and \frak switches to this
%
%%%  Add further fonts families here (using the same format) if you need
%    them.  The def of hexnumber is optional (it is only used for
%    \mathchardef 's).
%
%      Now we need to define the standard fonts (which are
%      already defined at 10,7 and 5 point) at 13,9 and 6 point:
%
%      Roman fonts:
\font\thirtnrm=cmr10 scaled 1315    %%%  ** see note above
\font\ninerm=cmr9                   \font\sixrm=cmr6   
%%%%%%  (add further sizes here if you need them)
%
%      Math italic fonts
\font\thirtni=cmmi10 scaled 1315    %%%  ** see note above 
\font\ninei=cmmi9                   \font\sixi=cmmi6  
%%%%%%  (add further sizes here if you need them)
%
%     Symbol fonts
\font\thirtnsy=cmsy10 scaled 1315   %%%  ** see note above
\font\ninesy=cmsy9                  \font\sixsy=cmsy6  
%%%%%%  (add further sizes here if you need them)
%
%     Bold face
\font\thirtnbf=cmbx10 scaled 1315   %%%  ** see note above 
\font\ninebf=cmbx9                  \font\sixbf=cmbx6  
%%%%%%  (add further sizes here if you need them)
%
%     The maths extension font (only defined at text size)
%
\font\thirtnex=cmex10 scaled 1315   %%%  ** see note above
\font\nineex=cmex9                  
%%%%%%  (add further sizes here if you need them)
%
%     Finally three fonts (text italic, slanted and typewriter type)
%     which are also only defined at text size
%
\font\thirtnit=cmti10 scaled 1315  %%%  ** see note above 
\font\nineit=cmti9                  
%%%%%%  (add further sizes here if you need them)
%
\font\thirtnsl=cmsl10 scaled 1315  %%%  ** see note above 
\font\ninesl=cmsl9                  
%%%%%%  (add further sizes here if you need them)
%
\font\thirtntt=cmtt10 scaled 1315  %%%  ** see note above 
\font\ninett=cmtt9                  
%%%%%%  (add further sizes here if you need them)
%
%
%     Now come the two main macros.  What  \small  does is to
%     change all the families of fonts from normal size which is
%     10,7,5  (ie 10pt text, 7pt subscript, 5pt subsubscript)
%     to 9,6,5.  \large  similarly changes to  13,9,7.  To make
%     other size changing macros, choose your three sizes, add
%     font size definitions if necessary and make the obvious changes
%     to one of these macros.  Change  \normalbaselineskip  and
%     \strutbox  dimensions to appropriate sizes as well.  To
%     add further fonts, insert them in each macro, using the
%     AMS fonts as a model.
%      
%
\def\small{%
%
%   redefine the sizes of the roman fonts :
%
\textfont0=\ninerm \scriptfont0=\sixrm \scriptscriptfont0=\fiverm
\def\rm{\fam0\ninerm}%       % ( \rm  sets \ninerm  in text mode
%                            %  and \fam0 in math mode)
%
%   and the math italic fonts :
%
\textfont1=\ninei \scriptfont1=\sixi \scriptscriptfont1=\fivei
%
%   and the symbol fonts :
%
\textfont2=\ninesy \scriptfont2=\sixsy \scriptscriptfont2=\fivesy
%
%   There is only one math extension font :
%
\textfont3=\nineex \scriptfont3=\nineex \scriptscriptfont3=\nineex
%
%   Next the bold font (named rather than numbered) :
%
\textfont\bffam=\ninebf \scriptfont\bffam=\sixbf
\scriptscriptfont\bffam=\fivebf \def\bf{\fam\bffam\ninebf}%
%
%   and the three text-only fonts : 
%
\textfont\itfam=\nineit \def\it{\fam\itfam\nineit}%
\textfont\slfam=\ninesl \def\sl{\fam\slfam\ninesl}%
\textfont\ttfam=\ninett \def\tt{\fam\ttfam\ninett}%
%
%   Now the three new families of AMS fonts :
%
%   AMS symbol-a
%
\textfont\msafam=\ninemsa \scriptfont\msafam=\sixmsa
\scriptscriptfont\msafam=\fivemsa \def\msa{\fam\msafam\ninemsa}%         
%
%   AMS symbol-b
%
\textfont\msbfam=\ninemsb \scriptfont\msbfam=\sixmsb
\scriptscriptfont\msbfam=\fivemsb \def\msb{\fam\msbfam\ninemsb}%         
%
%   Euler-Fraktal font
%
\textfont\eufmfam=\nineeufm  \scriptfont\eufmfam=\sixeufm
\scriptscriptfont\eufmfam=\fiveeufm \def\frak{\fam\eufmfam\nineeufm}%
%
%%%  Add further fonts families here if you need them.
%
%    Reset \normalbaselineskip and \strubox :
%
\normalbaselineskip=11pt%
\setbox\strutbox=\hbox{\vrule height8pt depth3pt width0pt}%
%
%    Set \normalbaselines and \rm (roman) as defaults :
%
\normalbaselines\rm
%
%    Reset some of the basic vertical skips:
%
\stdskip=4pt plus2pt minus2pt    
\medskipamount=\stdskip          
\parskip=\stdskip                
\abovedisplayskip=\stdskip       
\belowdisplayskip=\stdskip       
\def\ppar{\par\goodbreak\vskip 6pt plus 3pt minus 3pt}%     
%
%   And finally reset the size of section heads (see below):
%
\def\section##1{\global\advance\sectionnumber by 1
\vskip-\lastskip\penalty-800\vskip 20pt plus10pt minus5pt 
\egroup{\bf\number\sectionnumber\quad##1}\bgroup\small         
\vskip 6pt plus3pt minus3pt
\nobreak\resultnumber=1}%      % Reset resultnumber at start of section
}    %%%   End of  \small  macro      
%
%   Two useful abbreviations to keep track of \small material:
\def\beginsmall{\bgroup\small}
\let\endsmall\egroup
%
%
%    The \large  macro is similar (comments abbreviated):
%
%
\def\large{%
\textfont0=\thirtnrm \scriptfont0=\ninerm \scriptscriptfont0=\sevenrm
\def\rm{\fam0\thirtnrm}%
\textfont1=\thirtni \scriptfont1=\ninei \scriptscriptfont1=\seveni
\textfont2=\thirtnsy \scriptfont2=\ninesy \scriptscriptfont2=\sevensy
\textfont3=\thirtnex \scriptfont3=\thirtnex \scriptscriptfont3=\thirtnex
\textfont\bffam=\thirtnbf \scriptfont\bffam=\ninebf
\scriptscriptfont\bffam=\sevenbf \def\bf{\fam\bffam\thirtnbf}%
\textfont\itfam=\thirtnit \def\it{\fam\itfam\thirtnit}%
\textfont\slfam=\thirtnsl \def\sl{\fam\slfam\thirtnsl}%
\textfont\ttfam=\thirtntt \def\tt{\fam\ttfam\thirtntt}%
%   AMS symbol-a  :
\textfont\msafam=\thirtnmsa \scriptfont\msafam=\ninemsa
\scriptscriptfont\msafam=\sevenmsa \def\msa{\fam\msafam\thirtnmsa}%         
%   AMS symbol-b  :
\textfont\msbfam=\thirtnmsb \scriptfont\msbfam=\ninemsb
\scriptscriptfont\msbfam=\sevenmsb \def\msb{\fam\msbfam\thirtnmsb}%         
%   Euler-Fraktal font :
\textfont\eufmfam=\thirtneufm  \scriptfont\eufmfam=\nineeufm
\scriptscriptfont\eufmfam=\seveneufm \def\frak{\fam\eufmfam\teneufm}%
%%%% Add further fonts families here if you need them.
%   Reset \normalbaselineskip and \strubox and initialise :
\normalbaselineskip=16pt%
\setbox\strutbox=\hbox{\vrule height11.5pt depth4.5pt width0pt}%
\normalbaselines\rm}%
\let\Large\large   %  for compatibility with latex
%
%   The next two lines define commonly used switches for
%   blackboard bold (\Bbb) and gothic type (\goth).  The   
%   \Bbb  switch is set to work in the same way as in amstex
%   and switches only the next character to blackboard bold.
\def\Bbb#1{{\msb#1}}

%
%   To use the new AMS fonts you can either use the control
%   sequences \msa, \msb (alias \Bbb) and \frak (alias \goth) eg :

   % see the msam font table
%
%   or, more generally, make \mathchardef's (cf Knuth p155) eg :
\mathchardef\plussquare="0\hexa01
\mathchardef\nge="3\hexb0B
\mathchardef\maltesecross="0\hexa7A
\mathchardef\del="0\hexf01
%
%   or you can use the amstex names for all the new symbols by
%   inserting the line  \input amsnames  in your file directly
%   after \input gtmacros. 
%   This presupposes that you have collected a copy of the file
%   amsnames.tex  from the  gt/info/macros  ftp directory.
%
%
%   Finally we need a small capital font (for author(s)) :
%
\font\sc=cmcsc10
%
%%%%%%%%%%%%%%%%%       END OF FONT MACROS     %%%%%%%%%%%%%
%
%
%                 Knuth's \square macro :
%
\def\sqr#1#2{{\vcenter{\vbox{\hrule  height.#2truept
	\hbox{\vrule width.#2truept height#1truept 
	\kern#1truept \vrule width.#2truept}
	\hrule height.#2truept}}}}
%
    %   A small square for end-of-proofs. 
%                  %   (Define other size squares by varing the
%                  %   the two numbers.)
%
%
%      Style macros for section heads, theorem statements etc :
%   
%
\newcount\sectionnumber            %%%  Allocate registers to take
\newcount\resultnumber             %%%  section and result numbers.
\sectionnumber=0\resultnumber=1    %%%  Set these registers to 0 and 1
%
%   The \section macro produces a \large bold faced section heading
%   numbered to the left.  Pagebreaks are encouraged before the
%   start of the section and discouraged directly after the heading.
%   Typical use  \section{First steps}  with typical result :
%
%    1  First Steps     (set bold and \large)
%
\def\section#1{\global\advance\sectionnumber by 1
\xdef\nextkey{\number\sectionnumber}%      (used by the \key macro)
\vskip-\lastskip\penalty-800\vskip 20pt plus10pt minus5pt 
{\large\bf\number\sectionnumber\quad#1}         
\vskip 8pt plus4pt minus4pt
\nobreak\resultnumber=1}      % Reset resultnumber at start of section
%
%
%
%   Next a macro to set subheadings (like the  \section  macro
%   but without the number, with less space and set in standard size).
%
%   Typical use :  \sh{Example formats}
%
\def\sh#1{\vskip-\lastskip\ppar{\bf #1}\par\nobreak\medskip}         
%
%   The \proc ... \endproc macros ("proclaim") are for setting theorems, 
%   lemmas, conjectures etc with automatic numbering.  Typical use :    
%  
%    \proc{Theorem}Every lemon is yellow.\endproc
%
%   Typical result :
%     
%    Theorem 3.4  Every lemon is yellow.   

%   (with Theorem 3.4 set bold and a \stdspace of space before the 
%   statement set in slanted type).
%
\def\proc#1{\xdef\nextkey{\number\sectionnumber.\number\resultnumber}%
\vskip-\lastskip\ppar\bf%
\noindent#1\ \number\sectionnumber.\number\resultnumber
\stdspace\sl\global\advance\resultnumber by 1\ignorespaces}
 
%
%  The \prf ... \endprf macros are for setting proofs.  The code
%  for \prf includes the code for \endproc, so there is no need to
%  type \endproc if the theorem is followed immediatedly by a proof.
%
                            %  For start of proofs  
   %  For end (or absence) of proofs
                 %  extra vertical space)  
        %  For start of proof with alternative name
              %  \endproof is an alias for \endprf
%
%   Typical uses :    
%  
%    \proc{Theorem}Every lemon is yellow. \qed\endproc
%
%    \proc{Theorem}Every lemon is yellow.
%    \prf Use your eyes. \endprf
%
%    \proc{Theorem}Every lemon is yellow.
%    \proof{Proof of theorem} Use your eyes. \endprf
%
%   The next macro is a variant of the \proc macro.  It has
%   exactly the same result except that it omits the number.
%
%   Typical use :  
%    
%    \proclaim{Conjecture}Some oranges are yellow.\endproc
%

%
%   The next macro is a further variant for remarks, definitions etc.   
%   It omits the number and does not switch on slanted type.  
%  
%   Typical use :
%
%    \rk{Remark}Some lemons are thick-skinned.\endrk
%

%
%   The next macro is for numbering equations etc, \label  produces the 
%   correct number  x.y  and advances the resultnumber register
%
%   Typical use :
%
%     $$fx=7\eqno{\bf\label}$$
%
%   result :
%
%                           fx = 7                           3.5
%
\def\label{\xdef\nextkey{\number\sectionnumber.\number\resultnumber}%
\number\sectionnumber.\number\resultnumber
\global\advance\resultnumber by 1}
%
%
%
%   The next macros are to automate external references.  To use them 
%   type \reflist ..... \endreflist near the beginning of your paper, 
%   where  .... is the list of references in alphabetical order 
%   and in  the form  \key{KEY}  reference    where "KEY" is a 
%   string of characters which reminds you of the reference.   
%   Separate  references with a blank line or a \par.   Eg 
%
%     \reflist
%
%     ..... more references ....
%
%     \key{Kn-84} {\bf D Knuth}, {\it The TeXbook}, Addison--Wesley (1984)
%
%     ..... more references ....
%
%     \endreflist
%
%   Then type  \references  where you wish the references to be printed
%   (normally near the end of the paper).  To refer to Knuth type
%   for example    see Knuth [\ref{Kn-84}, page 320]   and the correct
%   numerical reference will be printed.  Edit the \references macro
%   to change the formatting (if desired).
%   There is an alternative \refkey for \key, provided your KEY contains
%   only letters.  The syntax is:
%
%     \reflist
%
%     ..... more references ....
%
%     \refkey\Knuth  {\bf D Knuth}, {\it The TeXbook}, Addison--Wesley (1984)
%
%     ..... more references ....
%
%     \endreflist
%
%   \key{Knuth}  has exactly the same maening as \refkey\Knuth and you
%   can mix the two syntaxes if you want.  But \refkey\Kn-84
%   would not work.  It would set Kn as the KEY and -84 would get printed!
%
\newcount\refnumber              %  Register for reference numbers
\refnumber=1                     %  set initially to 1.
\long\def\reflist#1\endreflist{%
\long\def\thereflist{#1}{\def\refkey##1##2\par{\xdef##1{\number\refnumber}%
\global\advance\refnumber by 1}%
\def\key##1##2\par{\expandafter\xdef%
\csname##1\endcsname{\number\refnumber}%
\global\advance\refnumber by 1}#1\par}}
\long\def\references{%
\penalty-800\vskip-\lastskip\vskip 15pt plus10pt minus5pt 
{\large\bf References}\ppar %`References' is set \large bold with space around.
{\leftskip=25pt\frenchspacing    % The list of references is set 
\small\parskip=3pt plus2pt       % \small  with small spaces between,
\def\refkey##1##2\par{\noindent  % numbers in [,]'s and set just to the
\llap{[##1]\stdspace}\ignorespaces##2\par}         % left of a 25pt margin.
\def\key##1##2\par{\noindent  
\llap{[\ref{##1}]\stdspace}\ignorespaces##2\par}  
\def\,{\thinspace}\thereflist\par}}
%
%   Next a footnote macro (with automatic numbering) which sets the
%   footnote  \small.
%
%   Typical use :
%         ..... are yellow.\fnote{By yellow here we mean Britsh
%    Standard colour BS3320.} 
%
\newcount\footnotenumber         % Register for footnote number
\footnotenumber=1                % set initially to 1
\def\fnote#1{\xdef\nextkey{\number\footnotenumber}%
{\small\ifnum\footnotenumber>9\parindent=14pt%
\else\parindent=10pt\fi\footnote{$^{\number\footnotenumber}$}%
{\hglue-5pt#1}\global\advance\footnotenumber by 1}}
%
%
%   Next macros for handling figures with automatic numbering (using 
%   TeX's \midinsert to float the figure to a suitable place).
%   
%   The \figure ... \endfigure macro centres the figure and adds
%   an automatically numbered label  Figure XX  after it.
%
%   If you have a caption, then type \caption{caption text} 
%   somewhere between \figure and \endfigure.  The macro
%   will then add  Figure XX: caption text  after the figure.
%
%   If you want an unnumbered or uncentred figure, then use TeX's raw 
%       \midinsert Figure instructions \endinsert  
%   and if you want a numbered figure label in the same style then
%   use \caption{caption text} outside of  \figure ... \endfigure.
%
%   If you need just the label Figure XX  outside of  \figure ... \endfigure
%   then type  \figurelabel .
%
\newcount\figurenumber          % register for figure number
\figurenumber=1                 % set initially to 1
\def\caption#1{\xdef\nextkey{\number\figurenumber}%
\cl{\small Figure \number\figurenumber: #1}%
\global\advance\figurenumber by 1}
\def\figurelabel{\xdef\nextkey{\number\figurenumber}%
\cl{\small Figure \number\figurenumber}%
\global\advance\figurenumber by 1}
\long\def\figure#1\endfigure{{\xdef\nextkey{\number\figurenumber}%
\let\captiontext\relax\def\caption##1{\xdef\captiontext{##1}}%
\midinsert\cl{\ignorespaces#1\unskip\unskip\unskip\unskip}\vglue6pt\cl{\small 
Figure \number\figurenumber\ifx\captiontext\relax\else: \captiontext
\fi}\endinsert\global\advance\figurenumber by 1}}
%
%   Macros for self-correcting internal references.
%
%   There are two macros  \key{KEY}  and  \ref{KEY} .
%
%   The \key macro sets up KEY as a key for whatever number is 
%   being referenced and the \ref macro converts the KEY into 
%   that number.  Type \key after a  \section or \proc or 
%   \label or \fnote or \figure or \caption or \figurelabel .
%
%   Example:
%
%       \section{Introduction}\key{intro}
%       \proc{Theorem}\key{MainTh}Lemons are yelloy\endproc
%       Here we follow\fnote{Follow in the sense of Dickens}
%       \key{Dickens-note}the crowd ....  
%
%       In section \ref{intro}
%       we stated theorem \ref{mainTh} and noted (see footnote 
%       \ref{Dickens-note}) ...
%
\def\nextkey{??}   %  initialise \nextkey (which is reset by all the
%                     numbering macros)
%
\def\key#1{\expandafter\xdef\csname #1\endcsname{\nextkey}}
\def\ref#1{\expandafter\ifx\csname #1\endcsname\relax
\immediate\write16{Reference {#1} undefined}??\else
\csname #1\endcsname\fi}
%
%   Note:  If the KEY contains only letters then \KEY has exactly the
%   same meaning as \ref{KEY} so in the example you could have:
%
%       In section \intro\ we ....
%
%   The \key will work at any time after the macro which sets the
%   number, provided no other macro which sets a number has been used. 
%
%   Macros for forward references:
%              =======
%   The \key \ref macros ONLY work for backwards references.  If you  
%   want to use forwards references, then type \useforwardrefs  near
%   the beginning of your file.  The KEY's are then stored in an
%   auxiliary  .ref  file and you then suffer the same disadvantage as
%   when using LaTeX that you must TeX the file twice to get
%   the references correct.
%
%   To use a forward ref type \ref{KEY}.  (You can type the
%   alternative  \KEY  but you'll get an error on first TeX'ing 
%   if the \KEY is not yet defined.) 
%
%   The macro also allows external references to be listed at the end 
%   of the file (if you wish to).  (Indeed they can be typed anywhere
%   before the \references command.)  You can combine the reference list
%   and the \references command by typing the references (using the
%   same syntax as before) between the commands \biblio and \endbiblio 
%   (don't type \references or they'll be printed twice).
%
\newread\gtinfile
\newwrite\gtreffile
\def\useforwardrefs{
\openin\gtinfile\jobname.ref
\ifeof\gtinfile
\closein\gtinfile
\immediate\write16{No file \jobname.ref}
\else
\closein\gtinfile
\input \jobname.ref
\fi
\immediate\openout\gtreffile \jobname.ref
%
%   Adapt \key :
%
\def\key##1{{\def\\{\noexpand}%
\expandafter\xdef\csname ##1\endcsname{\nextkey}%
\immediate\write\gtreffile{\\\expandafter\\\def\\\csname ##1\\\endcsname%
{\nextkey}}}}
%
%  Adapt macros for external references:  
%
\long\def\reflist##1\endreflist{%
\long\def\thereflist{##1}{\def\refkey####1####2\par{\xdef####1{%
\number\refnumber}{\def\\{\noexpand}\immediate\write\gtreffile
{\\\def\\####1{\number\refnumber}}}\global\advance\refnumber by 1}%
\def\key####1####2\par{\expandafter\xdef%
\csname####1\endcsname{\number\refnumber}%
{\def\\{\noexpand}\immediate\write\gtreffile
{\\\expandafter\\\def\\\csname ####1\\\endcsname{\number\refnumber}}}
\global\advance\refnumber by 1}##1\par}}
\long\def\biblio##1\endbiblio{\reflist##1\endreflist\references}%
%
%  Adapt obselete key macros (\numkey, \seckey and \figkey):
%
\def\numkey##1{{\def\\{\noexpand}%
\xdef##1{\number\sectionnumber.\number\resultnumber}
\immediate\write\gtreffile{\\\def\\##1%
{\number\sectionnumber.\number\resultnumber}}}}
\def\seckey##1{{\def\\{\noexpand}\xdef##1{\number\sectionnumber}
\immediate\write\gtreffile{\\\def\\##1{\number\sectionnumber}}}}
\def\figkey##1{\xdef##1{\number\figurenumber}%
{\def\\{\noexpand}\immediate\write\gtreffile%
{\\\def\\##1{\number\figurenumber}}}
\number\figurenumber\global\advance\figurenumber by 1}
}   %  end of \useforwardrefs
%
%
%   The next five macros are obselete and have been superseeded by
%   the general \key macro above.  They are included merely to 
%   maintain backward compatibility for the package:
%
%
\def\figkey#1{\xdef#1{\number\figurenumber}%
\number\figurenumber\global\advance\figurenumber by 1}
\def\fig#1#2\endfig{%
\midinsert\cl{#2}\vglue6pt\cl{\small Figure #1}\endinsert}
\def\newfig{\number\figurenumber\global\advance\figurenumber by 1}
\def\numkey#1{\xdef#1{\number\sectionnumber.\number\resultnumber}}
\def\seckey#1{\xdef#1{\number\sectionnumber}}
%
%   End of obselete macros.
%
%
%   The next macro is a version of the verbatim macro given by Knuth.
%
%   This macro produces a "verbatim" printout of
%   any ASCII string which does not contain the symbol "
%   (TeX files do not usually contain " 's).
%   More precisely, everything between consecutive pairs
%   of " 's is printed verbatim in the typewriter font cmtt.
%   For an explanation of how the macro works, see Knuth pp 420-1.
%
%   There are two switches: \verb (which switches the macro on)
%   and \brev which switches the macro off (the default).  When
%   the macro is switched off the symbol " has its usual 
%   meaning for TeX.  To use the macro, type \verb before use
%   and the use " to switch verbatim on and off.  Be careful
%   not to use " for any other purpose.  There is no need to
%   switch the macro off again unless you need to use " for
%   some other purpose (eg making  \mathchardef 's).  Note 
%   that the macro MUST BE OFF before inputting  amsnames.tex .
%
%   Whether the macro is on or off you can always use the
%   control sequence \dq (double quote) for " e.g.
%   \mathchardef\sum=\dq1350  is perfectly valid.
%   The control sequence \ttq is an abbreviation for
%   {\tt\dq}.  Thus "\ttq" will produce " (in cmtt)
%   inside a verbatim quote.
%
%
   %  define a code for " so it can be used when \verb is on
  %  code for " in cmtt
%
\def\verb{\catcode`\"=\active}       %  The main
\def\brev{\catcode`\"=12}            %  switches.
\brev                                %  Prime switches and
\verb                                %  switch on.
{\obeyspaces\gdef {\ }}              
{\catcode`\`=\active\gdef`{\relax\lq}}
\def"{%
\begingroup\baselineskip=12pt\def\par{\leavevmode\endgraf}%
\tt\obeylines\obeyspaces\parskip=0pt\parindent=0pt%
\catcode`\$=12\catcode`\&=12\catcode`\^=12\catcode`\#=12%
\catcode`\_=12\catcode`\~=12%
\catcode`\{=12\catcode`\}=12\catcode`\%=12\catcode`\\=12%
\catcode`\`=\active\let"\endgroup}
\brev      %   Finally switch the macro off (for safety)
%
%   Macros for itemised lists.   Typical use :
%    
%    \items
%    \item{(i)}Colours must be defined.
%    \item{(ii)}Colour cards may not be cited.
%    \enditems
%
%   Result :
%
%    (i)  Colours must be defined. 
%   (ii)  Colour cards may not be cited.
%
%
           % Start of itemised list         
         % end of itemised list   
\def\item#1{\par\leavevmode\llap{#1\stdspace}%
\ignorespaces}                             % labelled item
               % bulleted item.
%
%   The \quote ... \endquote macros are for typesetting quotations :
%

%
%   A few useful abbreviations :
%
    %  Colon with correct spacing for maps.
\def\np{\vfil\eject}         %  Forced page break (new page).
\def\nl{\hfil\break}         %  New line.
\def\cl{\centerline}         %  Centerline
        %  The journal title in recommended style
\def\gtm{{\mathsurround=0pt\it $\cal G\mskip-2mu$eometry \&\ 
$\cal T\!\!$opology $\cal M\mskip-1mu$onographs}}    %  for monographs
\def\agt{{\mathsurround=0pt\it$\cal A\mskip-.7mu$lgebraic \&\ 
$\cal G\mskip-2mu$eometric $\cal T\!\!$opology}}  % AGT
%
%    Finally some macros for automatic title page or header generation.
%    To use them type your header information using the following  
%    example as a guide :
%
%    Note that \\ is used as standard separator (for lines in \title and
%    \address, between authors and between email addresses or URL's)
%    and that \email, \url and \secondaddress are optional.
%

% Example:  \title{A short spoof paper\\with a two-line title}
% =======   \authors{Albert Einstein\\Leonardo da Vinci}
%           \address{IAS\\Princeton}\secondaddress{Renaissance\\Venice}
%           \email{ae@ias.princeton.edu\\ldv@ren.ven.hist}
%           \abstract 
%           A short spoof paper with a very short abstract.
%           \endabstract 
%           \primaryclass{00-01, 00-02}\secondaryclass{68-00, 68-01}
%           \keywords{Short, spoof, paper}
%           \maketitlepage
%
%
%    The title page or header will then be generated automatically.
%
%
%    Define the various ingredients of the title page:
%
\def\title#1{\def\thetitle{#1}}

\def\author#1{\edef\previousauthors{\theauthors}
 \ifx\theauthors\relax\def\theauthors{#1}\else
 \def\theauthors{\previousauthors\par#1}\fi}

\let\authors\author        % aliases
\def\address#1{\edef\previousaddresses{\theaddress}
 \ifx\theaddress\relax\def\theaddress{#1}\else
 \def\theaddress{\previousaddresses\par\vskip 2pt\par#1}\fi}
                             % alias
\def\secondaddress#1{\edef\previousaddresses{\theaddress}
 \ifx\theaddress\relax\def\theaddress{#1}\else
 \def\theaddress{\previousaddresses\par{\rm and}\par#1}\fi}   

\def\email#1{\edef\previousemails{\theemail}
 \ifx\theemail\relax\def\theemail{#1}\else
 \def\theemail{\previousemails\hskip 0.75em\relax#1}\fi}
  % aliases
\def\secondemail#1{\edef\previousemails{\theemail}
 \ifx\theemail\relax\def\theemail{#1}\else
 \def\theemail{\previousemails\hskip 0.75em{\rm and}\hskip 0.75em
 \relax#1}\fi}
\def\url#1{\edef\previousurls{\theurl}
 \ifx\theurl\relax\def\theurl{#1}\else
 \def\theurl{\previousurls\hskip 0.75em\relax#1}\fi}
      % aliases
\def\secondurl#1{\edef\previousurls{\theurl}
 \ifx\theurl\relax\def\theurl{#1}\else
 \def\theurl{\previousurls\hskip 0.75em{\rm and}\hskip 0.75em
 \relax#1}\fi}
\long\def\abstract#1\endabstract{\long\def\theabstract{#1}}
\def\primaryclass#1{\def\theprimaryclass{#1}}
                        % alias
\def\secondaryclass#1{\def\thesecondaryclass{#1}}
\def\keywords#1{\def\thekeywords{#1}}
%
%  Set \\ to \par and title page items to \relax to initialise macros :
%
\let\\\par\let\thetitle\relax\let\theshorttitle\relax
\let\theauthors\relax\let\theshortauthors\relax
\let\theaddress\relax\let\theshortaddress\relax
\let\theemail\relax\let\theurl\relax
\let\theabstract\relax\let\theprimaryclass\relax
\let\thesecondaryclass\relax\let\thekeywords\relax
%
%
%
%   Basic title page layout (edit this macro if you
%   wish to adjust the title page layout) :
%
\long\def\maketitlepage{    % start of definition of \maketitlepage

\vglue 0.2truein   % top margin

% title :
%
{\parskip=0pt\leftskip 0pt plus 1fil\def\\{\par\smallskip}{\large
\bf\thetitle}\par\medskip}   

\vglue 0.15truein 

% authors :
%
{\parskip=0pt\leftskip 0pt plus 1fil\def\\{\par}{\sc\theauthors}
\par\medskip}%
 
\vglue 0.1truein 

% address(es) email's and URL's (with switches to detect whether the
% optional items have been used) :
%
{\small\parskip=0pt
{\leftskip 0pt plus 1fil\def\\{\par}{\sl\theaddress}\par}
\ifx\theemail\relax\else  % email address?
\vglue 5pt \def\\{\stdspace{\rm and}\stdspace} 
\cl{Email:\stdspace\tt\theemail}\fi
\ifx\theurl\relax\else    % URL given?
\vglue 5pt \def\\{\stdspace{\rm and}\stdspace} 
\cl{URL:\stdspace\tt\theurl}\fi\par}

\vglue 7pt 

{\bf Abstract}

\vglue 5pt

\theabstract

\vglue 7pt 

{\bf AMS Classification numbers}\quad Primary:\quad \theprimaryclass\par

Secondary:\quad \thesecondaryclass

\vglue 5pt 

{\bf Keywords:}\quad \thekeywords

\np  % page break at the end of the title page

}    % end of definition of \maketitlepage
%
%    % \makeshorttitle (for general preprints) doesn't take a new page
%
\long\def\makeshorttitle{    % start of definition of \makeshorttitle

%\vglue 0.2truein   % top margin

% title :
%
{\parskip=0pt\leftskip 0pt plus 1fil\def\\{\par\smallskip}{\large
\bf\thetitle}\par\medskip}   

\vglue 0.05truein 

% authors :
%
{\parskip=0pt\leftskip 0pt plus 1fil\def\\{\par}{\sc\theauthors}
\par\medskip}%
 
\vglue 0.03truein 

% address(es) email's and URL's (with switches to detect whether the
% optional items have been used) :
%
{\small\parskip=0pt
{\leftskip 0pt plus 1fil\def\\{\par}{\sl\ifx\theshortaddress\relax
\theaddress\else\theshortaddress\fi}\par}
\ifx\theemail\relax\else  % email address?
\vglue 5pt \def\\{\stdspace{\rm and}\stdspace} 
\cl{Email:\stdspace\tt\theemail}\fi
\ifx\theurl\relax\else    % URL given?
\vglue 5pt \def\\{\stdspace{\rm and}\stdspace} 
\cl{URL:\stdspace\tt\theurl}\fi\par}

\vglue 10pt 

% abstract and classification numbers (with switches):

{\small\leftskip 25pt\rightskip 25pt{\bf Abstract}\stdspace\theabstract

{\bf AMS Classification}\stdspace\theprimaryclass
\ifx\thesecondaryclass\relax\else; \thesecondaryclass\fi\par
{\bf Keywords}\stdspace \thekeywords\par}
\vglue 7pt
}    % end of definition of \makeshorttitle
\let\maketitle\makeshorttitle        %% alias
%
%    %%%% \makeagttitle (for AGT) similar to \makeshorttitle but
%         with addresses omitted (they go at the end)
%
%%%% publication info and test defaults:

\def\volumenumber#1{\def\thevolumenumber{#1}}
\def\volumeyear#1{\def\thevolumeyear{#1}}
\def\pagenumbers#1#2{\def\startpage{#1}\def\finishpage{#2}}
\def\published#1{\def\publishdate{#1}}
\def\received#1{\def\receiveddate{#1}}
\def\revised#1{\def\reviseddate{#1}}
\let\reviseddate\relax
%% Defaults for authors to use to check layout
\volumenumber{X}
\volumeyear{20XX}
\pagenumbers{1}{XXX}
\published{XX Xxxember 20XX}

\long\def\makeagttitle{   %%% start of definition of \makeagttitle
\agt\hfill      %   Journal title (top left) 
%   logo placeholder (top right)
\hbox to 60truept{\vbox to 0pt{\vglue -14truept{\bf [Logo here]}\vss}\hss}
\break
{\small Volume \thevolumenumber\ (\thevolumeyear)
\startpage--\finishpage\nl
Published: \publishdate}

\vglue .2truein

% title
{\parskip=0pt\leftskip 0pt plus 1fil\def\\{\par\smallskip}{\large
\bf\thetitle}\par\medskip}   
\vglue 0.05truein 

% authors :
%
{\parskip=0pt\leftskip 0pt plus 1fil\def\\{\par}{\sc\theauthors}
\par\medskip}%
 
\vglue 0.03truein 

%  abstract and classification numbers:

{\small\leftskip 25truept\rightskip 25truept{\bf Abstract}\stdspace\theabstract

{\bf AMS Classification}\stdspace\theprimaryclass
\ifx\thesecondaryclass\relax\else; \thesecondaryclass\fi\par
{\bf Keywords}\stdspace \thekeywords\par}\vglue 7truept

}   %%%% end of definition of \makeagttitle

%%%%% Macro to typeset addresses (typically at the end of the paper)

\def\Addresses{\bigskip
{\small \parskip 0pt \leftskip 0pt \rightskip 0pt plus 1fil \def\\{\par}
\sl\theaddress\par\medskip \rm Email:\stdspace\tt\theemail\par
\ifx\theurl\relax\else\smallskip \rm URL:\stdspace\tt\theurl\par\fi}}

\def\agtart{%   Full mock-up of AGT article style (for authors to test with)
%  get print centerpage:
\hoffset 14truemm
\voffset 31truemm
\font\phead=cmsl9 scaled 950
\font\pnum=cmbx10 scaled 913
\font\pfoot=cmsl9 scaled 950
%  headline and footline
\headline{\vbox to 0pt{\vskip -4.5mm\line{\small\phead\ifnum
\count0=\startpage ISSN numbers are printed here
\hfill {\pnum\folio}\else\ifodd\count0\def\\{ }% 
\ifx\theshorttitle\relax\thetitle\else\theshorttitle\fi\hfill{\pnum\folio}
\else\def\\{ and }{\pnum\folio}\hfill\ifx\theshortauthors\relax\theauthors
\else\theshortauthors\fi\fi\fi}\vss}}
\footline{\vbox to 0pt{\vglue 0mm\line{\small\pfoot\ifnum\count0=\startpage
Copyright declaration is printed here\hfill\else
\agt, Volume \thevolumenumber\ (\thevolumeyear)\hfill\fi}\vss}}
%  force \agttitle
\let\maketitle\makeagttitle\let\makeshorttitle\makeagttitle}

%%%
%%%  This is gtmonout.tex. It contains routines which automatically
%%%  generate:  
%%%  (1)  the title page in correct format
%%%  (2)  TeX file for the reprint cover (cover.tex)
%%%  (3)  a short batch file to rename the reprint cover file (covname.bat)
%%%  (4)  a data file (paper.dat) for the final publication routine
%%%  (5)  a text version (paper.txt) for email to subscribers and abstract
%%%  (6)  an html version (paper.htm) for the html abstract.
%%%
%%%  Used with makemon.bat, this completely automates the final
%%%  production of GTM articles and reprint cover, and largely automates
%%%  the final publication. 
%%%
%%%  For instructions see gtoutput.txt.
%%%                                               Colin Rourke  26.7.02

%  test for latex or plain tex
\def\ifplaintex{\expandafter\ifx\csname documentclass\endcsname\relax}

%  get print centerpage:

\ifplaintex 
\hoffset 14truemm
\voffset 31truemm
\else
\headsep 23pt
\footskip 35pt
\hoffset -4truemm
\voffset 12.5truemm
\fi

        %  journal title in recommended style

\def\gtm{{\mathsurround=0pt\it $\cal G\mskip-2mu$eometry \&\ 
$\cal T\!\!$opology $\cal M\mskip-1mu$onographs}}    %  for monographs

\def\gtp{{\mathsurround=0pt\it $\cal G\mskip-2mu$eometry \&\ 
$\cal T\!\!$opology $\cal P\!$ublications}}  % GT publications

\def\recd{{\small Received:\qua\receiveddate\ifx\reviseddate\relax
\else\qquad Revised:\qua\reviseddate\fi\par}} 

%  define the various new ingredients of the title page and the data
%  output files

\def\volumenumber#1{\def\thevolumenumber{#1}}
\def\volumeyear#1{\def\thevolumeyear{#1}}
\def\volumename#1{\def\thevolumename{#1}}
\def\papernumber#1{\def\thepapernumber{#1}}
\def\pagenumbers#1#2{\def\startpage{#1}\def\finishpage{#2}}
\def\published#1{\def\publishdate{#1}}
\def\received#1{\def\receiveddate{#1}}
\def\revised#1{\def\reviseddate{#1}}
\def\accepted#1{\def\accepteddate{#1}}

%  initialise

\let\\\par
\let\thevolumenumber\relax\let\thepapernumber\relax
\let\thevolumeyear\relax\let\startpage\relax
\let\finishpage\relax\let\publishdate\relax\let\receiveddate\relax
\let\reviseddate\relax\let\accepteddate\relax\let\theasciititle\relax
\let\theasciiauthors\relax
\let\theasciiabstract\relax

\let\theerratum\relax\let\theasciiemail\relax
\let\theshortauthors\relax\let\theshorttitle\relax

%%% Define a few things for test purposes
\def\startpage{1}\def\finishpage{15}\def\thepapernumber{77}

%%%% Edit the next three lines for later volumes (or include definitions
%%%% with papers):
\volumenumber{2}
\volumename{Proceedings of the Kirbyfest}
\volumeyear{1999}

\long\def\maketitlep{   % start of definition of \maketitlep

\count0=\startpage

\gtm\nl        %   GT mongraphs (top left) 
{\small Volume \thevolumenumber: \thevolumename\nl 
\ifx\theerratum\relax\else Erratum \erratumnumber\nl\fi
Pages \startpage--\finishpage\nl}

\vglue 0.1truein   % top margin

% title
{\parskip=0pt\leftskip 0pt plus 1fil\def\\{\par\smallskip}{\ifplaintex\large
\else\Large\fi\bf\thetitle}\par\medskip}   
\vglue 0.05truein 

% authors :
%
{\parskip=0pt\leftskip 0pt plus 1fil\def\\{\par}{\sc\theauthors}
\par\medskip}%
 
\vglue 0.03truein 

%  abstract and classification numbers:

{\small\leftskip 25pt\rightskip 25pt{\bf Abstract}\stdspace\theabstract

{\bf AMS Classification}\stdspace\theprimaryclass
\ifx\thesecondaryclass\relax\else; \thesecondaryclass\fi\par
{\bf Keywords}\stdspace \thekeywords\par}\vglue 7pt

}   % end of definition of \maketitlep

%%% Headers and footers

\font\phead=cmsl9 scaled 950
\font\lhead=cmsl9 scaled 1050
\font\pnum=cmbx10 scaled 913
\font\lnum=cmbx10 
\font\pfoot=cmsl9 scaled 950
\font\lfoot=cmsl9 scaled 1050
\ifplaintex
\headline{\vbox to 0pt{\vskip -4.5mm\line{\small\phead\ifnum
\count0=\startpage ISSN 1464-8997 (on line)
1464-8989 (printed) \hfill {\pnum\folio}\else\ifodd\count0\def\\{ }% 
\ifx\theshorttitle\relax\thetitle\else\theshorttitle\fi\hfill{\pnum\folio}
\else\def\\{ and }{\pnum\folio}\hfill\ifx\theshortauthors\relax\theauthors
\else\theshortauthors\fi\fi\fi}\vss}}
\footline{\vbox to 0pt{\vglue 0mm\line{\small\pfoot\ifnum\count0=\startpage
Published \publishdate:\qua\copyright\ \gtp\hfill\else
\gtm, Volume \thevolumenumber\ (\thevolumeyear)\hfill\fi}\vss
}}
\else
\makeatletter
\def\@oddhead{{\small\lhead\ifnum\count0=\startpage ISSN 1464-8997 (on line)
1464-8989 (printed) \hfill {\lnum\number\count0}\else\ifodd\count0
\def\\{ }\ifx\theshorttitle\relax \thetitle \else\theshorttitle\fi\hfill
{\lnum\number\count0}\else\def\\{ and }{\lnum\number\count0}
\hfill\ifx\theshortauthors\relax 
\theauthors\else\theshortauthors\fi\fi\fi}}\def\@evenhead{@oddhead}
\def\@oddfoot{\small\lfoot\ifnum\count0=\startpage Published \publishdate:\qua\copyright\ \gtp\hfill\else
\gtm, Volume \thevolumenumber\ (\thevolumeyear)\hfill\fi}
\def\@evenfoot{@oddfoot}
\makeatother
\fi

\let\maketitlepage\maketitlep
\let\makeshorttitle\maketitlepage
\let\maketitle\maketitlepage

%\endinput  %%%comment out to create xxx header file

\newwrite\gtoutfile
\long\gdef\makeheadfile{  %%% start of definition of \makeheadfile
{\def\\{, }\def\s{ }
\immediate\openout\gtoutfile head.xxx
\immediate\write\gtoutfile{To: math@arxiv.org}
\immediate\write\gtoutfile{Subject: put OR rep NNNNN:ppppp}
\immediate\write\gtoutfile{--text follows this line--}
\immediate\write\gtoutfile{Proxy-for: \ifx\theasciiauthors\relax
\theauthors\else\theasciiauthors\fi\s<\ifx\theasciiemail\relax\theemail\else\theasciiemail\fi>}
\immediate\write\gtoutfile{\noexpand\\}
\immediate\write\gtoutfile{Authors: \ifx\theasciiauthors\relax
\theauthors\else\theasciiauthors\fi}
{\def\\{ }\immediate\write\gtoutfile{Title: \ifx\theasciititle\relax
\thetitle\else\theasciititle\fi}}
\immediate\write\gtoutfile{Subj-class: GT or SG, GR etc}
\immediate\write\gtoutfile{MSC-class: \theprimaryclass\ifx\thesecondaryclass\relax\else, \thesecondaryclass\fi}
\immediate\write\gtoutfile{Journal-ref: Geom. Topol. Monogr. \thevolumenumber\s
(\thevolumeyear) \startpage-\finishpage}
\immediate\write\gtoutfile{Comments: Published by Geometry and Topology Monographs at}
\immediate\write\gtoutfile{\s\s\s  http://www.maths.warwick.ac.uk/gt/GTMon\thevolumenumber/paper\thepapernumber.abs.html}
\immediate\write\gtoutfile{\noexpand\\}
\immediate\write\gtoutfile{}
\ifx\theasciiabstract\relax
\immediate\write\gtoutfile{\theabstract}\else
\immediate\write\gtoutfile{\theasciiabstract}\fi
\immediate\write\gtoutfile{}
\immediate\write\gtoutfile{\noexpand\\}
\immediate\write\gtoutfile{}
\immediate\closeout\gtoutfile}}  %%% end of definition of \makeheadfile

\def\maketitlepage{\maketitlep\makeheadfile}
\let\makeshorttitle\maketitlepage
\let\maketitle\maketitlepage

\volumenumber{4}
\volumename{Invariants of knots and 3-manifolds (Kyoto 2001)}
\volumeyear{2002}
\papernumber{20}
\pagenumbers{303}{311}
\received{5 December 2001}
\revised{19 February 2002}
\accepted{22 July 2002}
\published{13 October 2002}

\input epsf

\def\head #1.#2\endhead{\section{#2}}
\def\subhead#1\endsubhead{\sh{#1}}
\def\subsubhead#1\endsubsubhead{\sh{#1}}
\let\Cal\cal
\def\align#1\endalign{\eqalign{#1}}
\def\text#1{\hbox{#1}}

\reflist

\key{1} {\bf H~Murakami},
{\it Optimistic calculations about the Witten-Reshetikhin-Turaev invariants of
closed three-manifolds obtained from the figure-eight knot by integral Dehn surgeries},
\lq\lq Recent Progress Toward the Volume Conjecture", RIMS Kokyuroku 1172 (2000) 70--79

\key{2} {\bf H~Murakami}, {\bf J~Murakami}, {\bf M~Okamoto}, {\bf T~Takata}, {\bf Y~Yokota},
{\it Kashaev's conjecture and the Chern-Simons invariants of knots and links}, preprint

\key{3} {\bf W\,D~Neumann}, {\bf D~Zagier},
{\it Volumes of hyperbolic 3-manifolds}, Topology 24 (1985) 307--332

\key{4} {\bf Y~Yokota},
{\it On the volume conjecture for hyperbolic knots},
preprint available at {\tt http://www.comp.metro-u.ac.jp/\~{}jojo/volume-conjecture.ps}

\key{5} {\bf T~Yoshida},
{\it The $\eta$-invariant of hyperbolic 3-manifolds},
Invent. Math. 81 (1985) 473--514

\endreflist

\title{On the potential functions for the hyperbolic\\structures 
of a knot complement}                    
\authors{Yoshiyuki Yokota}                  
\email{jojo@math.metro-u.ac.jp}

\address{Department of Mathematics, Tokyo Metropolitan 
University\\Tokyo, 192-0397, Japan}                  

\abstract  
We explain how to construct certain potential functions
for the hyperbolic structures of a knot complement,
which are closely related to the analytic functions
on the deformation space of hyperbolic structures.
\endabstract

\primaryclass{57M50}                
\secondaryclass{57M25, 57M27}              
\keywords{Potential function, hyperbolicity equations, volume, 
Chern-Simons invariant}                    
\maketitle  

\cl{\small\it Dedicated to Professor Mitsuyoshi Kato for his 60'th birthday}

\head 1. Introduction \endhead

Let $M$ be the complement of a hyperbolic knot $K$ in $S^3$.
Through the study of {\it Kashaev's conjecture},
we have found a complex function
which gives the {\it volume} and the {\it Chern-Simons invariant}
of the complete hyperbolic structure of $M$ at the critical point
corresponding to the promised solution to the hyperbolicity equations for $M$,
see [2, 4] for details.

The purpose of this article is
to explain how to construct such complex functions
for the non-complete hyperbolic structures of $M$.
Such functions are closely related to
the analytic functions on the {deformation space} of the hyperbolic structures of $M$,
parametrized by the eigenvalue of the holonomy representation of the meridian of $K$,
which reveal a complex-analytic relation
between the volumes and the Chern-Simons invariants of the hyperbolic structures of $M$,
see [3, 5] for details.

In this note, we suppose $K$ is $5_2$ for simplicity
which is represented by the diagram $D$ depicted in Figure 1.

\figure
\epsfxsize0.7in\epsfbox{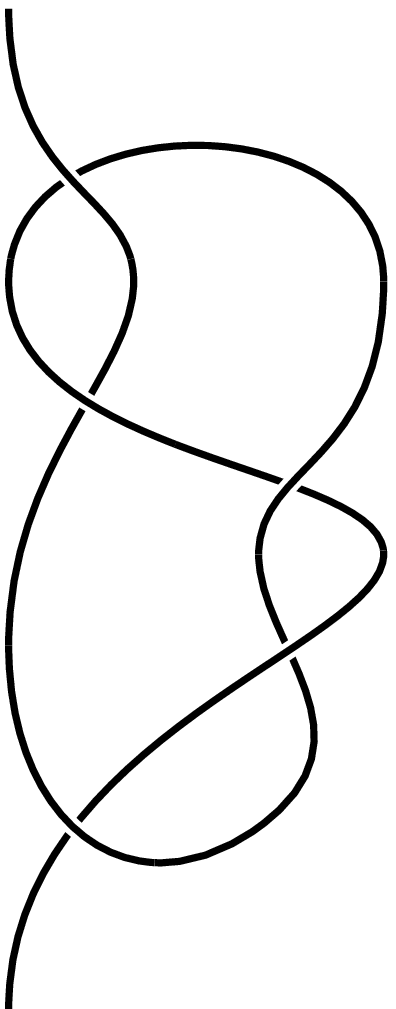}
\endfigure

\head 2. Geometry of a knot complement \endhead

\subsubhead 2.1 Ideal triangulations\endsubsubhead
We first review an ideal triangulation of $M$ due to D.~Thurston.
Let $\dot M$ denote $M$ with two poles $\pm\infty$ of $S^3$ removed.
Then, $\dot M$ decomposes into 5 ideal octahedra corresponding to the 5 crossings of $D$,
each of which further decomposes into 4 ideal tetrahedra around an axis, as shown in Figure 2.

\figure
\epsfxsize1in\epsfbox{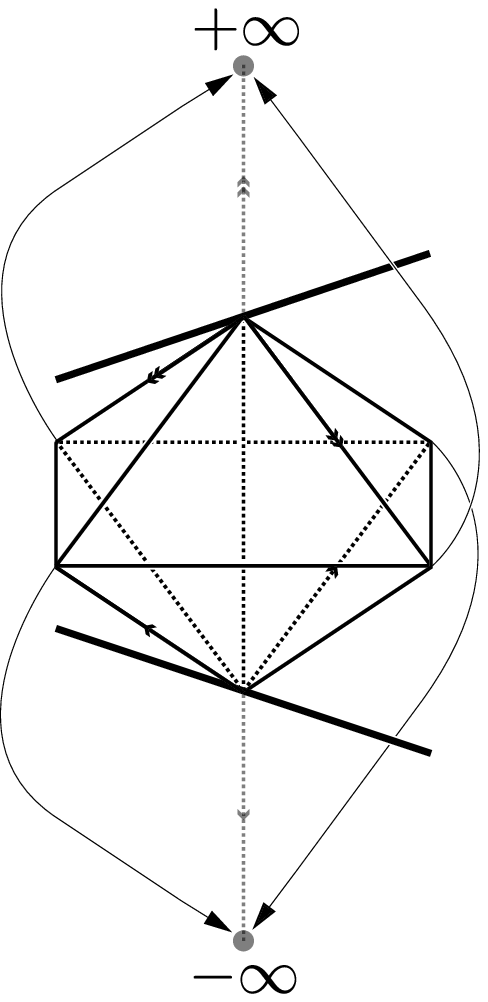}
\endfigure

In fact, we can reocover $\dot M$ by glueing adjacent tetrahedra as shown in Figure 3.

\figure
\epsfxsize2.5in\epsfbox{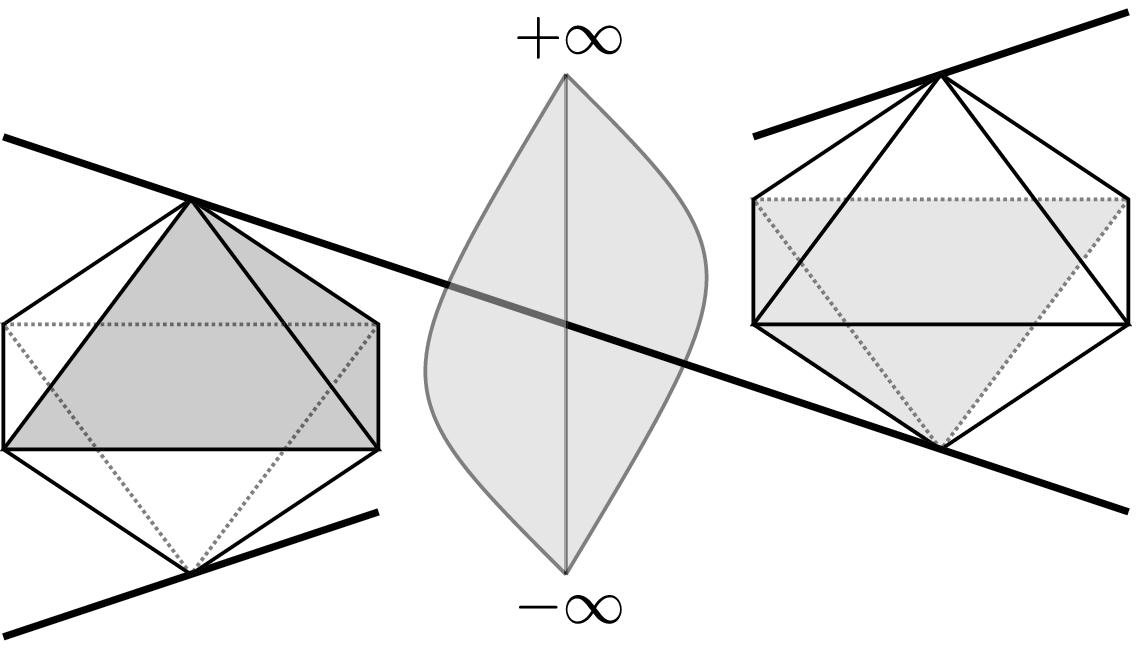}
\endfigure

As usual, we put a hyperbolic structure on each tetrahedron
by assigning a complex number, called {\it modulus},
to the edge corresponding to the axis as shown in Figure 4.
In what follows, we denote the tetrahedron with modulus $z$ by $T(z)$.

\figure
\epsfxsize0.9in\epsfbox{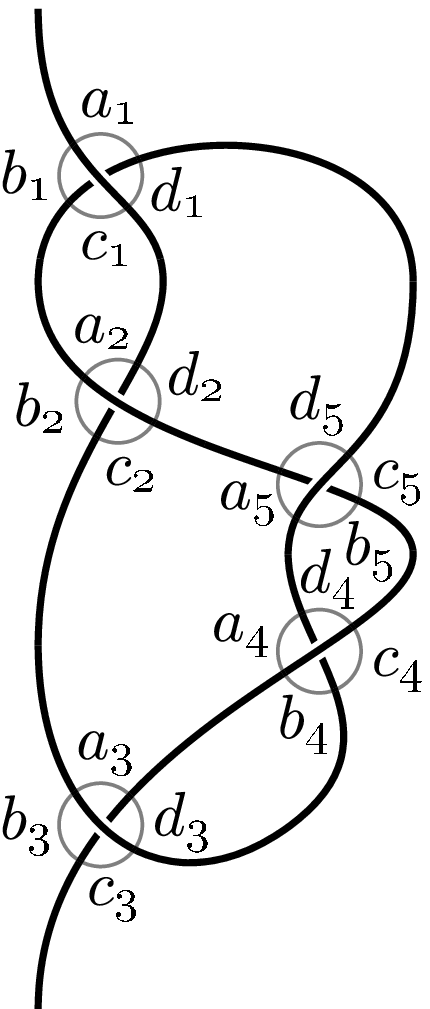}
\endfigure

Let $B$ be the intersection between $T(a_1)\cup T(b_1)$ and $T(b_3)\cup T(c_3)$.
Then, each of
$$
T(a_1),T(b_1),T(b_3),T(c_3)
$$
intersects $\partial N(B\cup K)$ in two triangles,
and they are essentially one-dimensional objects in $S^3\setminus N(B\cup K)$.
On the other hand, each of
$$
T(c_1),T(d_1),T(a_2),T(b_2),T(d_2),T(a_3),T(d_3),T(a_4),T(b_4),T(c_4),T(c_5)
$$
intersects $\partial N(B\cup K)$ in two triangles and one quadrangle,
and they are essentially two-dimensional objects in $S^3\setminus N(B\cup K)$.
Thus, by contracting these 15 tetrahedra,
we obtain an ideal triangulation $\Cal S$ of $M$ with
$$
T(c_2),T(d_4),T(a_5),T(b_5),T(d_5).
$$
Figure 5 exhibits the triangulation of $\partial N(B\cup K)$ induced by $\Cal S$,
where each couple of edges labeled with the same number are identified.

\figure
\epsfxsize4.5in\epsfbox{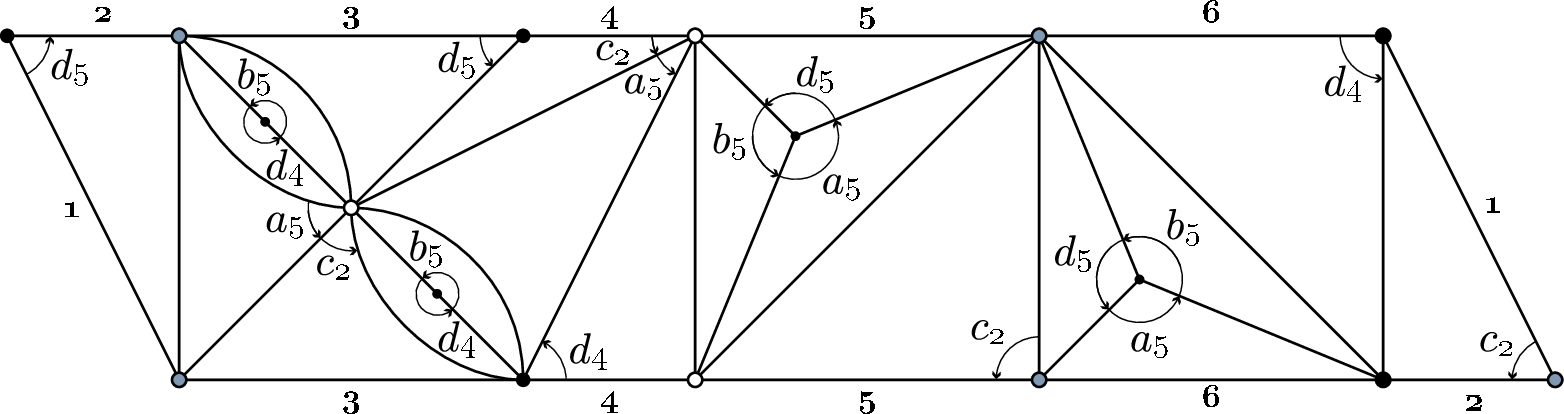}
\endfigure

\subsubhead 2.2 Hyperbolicity equations\endsubsubhead
If $c_2,d_4,a_5,b_5,d_5$ above give a hyperbolic structure of $M$,
the product of the moduli around each edge in $\Cal S$ should be 1,
which is called the {\it hyperbolicity equations}
and can be read from Figure 5 as follows.
$$
\align
&d_4b_5=\ a_5b_5d_5=1,
\cr
&{c_2a_5(1-1/d_4)\over1-d_4}\cdot
{(1-1/d_5)(1-1/c_2)(1-1/b_5)\over(1-a_5)(1-b_5)}=1,
\cr
&{c_2(1-1/a_5)\over(1-d_5)(1-c_2)}\cdot
{(1-1/d_5)(1-1/b_5)\over(1-d_5)(1-a_5)(1-d_4)}=1,
\cr
&{d_4(1-1/a_5)(1-1/d_4)\over1-b_5}\cdot{d_5(1-1/c_2)\over1-c_2}=1.
\endalign
$$
It is easy to observe that these equations are generated by
$$
d_4b_5=a_5b_5d_5=1,\ 
{c_2a_5\over d_4}={1-1/a_5\over(1-1/c_2)(1-d_5)}={(1-1/a_5)(1-d_4)\over1-b_5}={c_2\over d_5},
$$
which suggests to put
$$
c_2=y\xi,\ d_4=x/\xi,\ a_5=x/y,\ b_5=\xi/x,\ d_5=y/\xi
$$
and to rewrite the hyperbolicity equations as follows.
$$
{(1-y/x)(1-x/\xi)\over1-\xi/x}={1-y/x\over(1-y/\xi)(1-1/y\xi)}=\xi^2.
$$
Note that the variables $x,y$ correspond to the interior edges of a graph depicted in Figure 6,
which is $D$ with some edges deleted.

\figure
\epsfxsize 1in\epsfbox{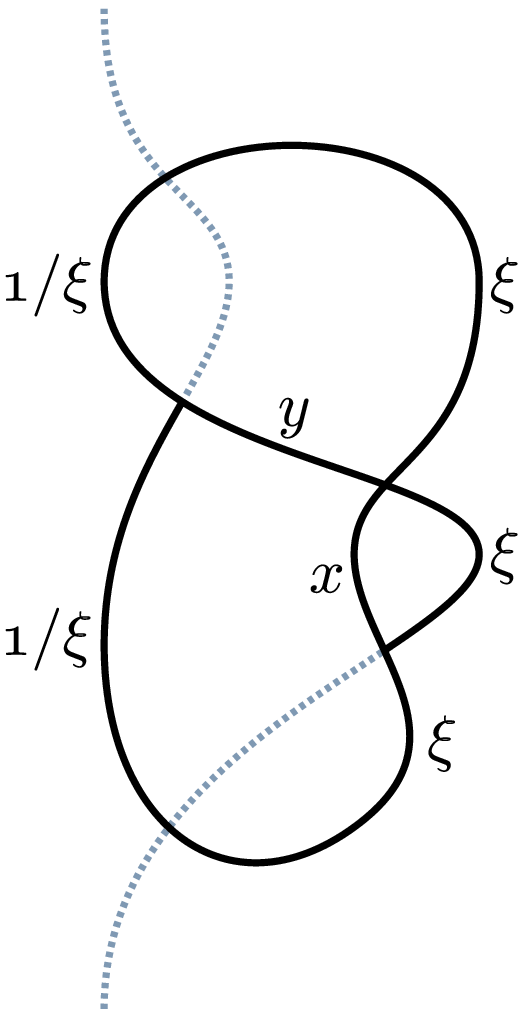}
\endfigure

A solution to the equations above determines a hyperbolic structure of $M$,
where $\xi$ is nothing but the eigenvalue of the holonomy representation of the meridian of $K$.
The set $\Cal D$ of such solutions is called
the {\it deformation space} of the hyperbolic structures of $M$
and can be parametrized by $\xi$
or the eigenvalue $\eta$ of the holonomy representation of the longitude of $K$.
In our example, $\eta$ is given by
$$
\eta
={y\xi^6\over x}\cdot(1-1/y\xi)
={y\xi^6\over x}\cdot{1-\xi/x\over(1-x/\xi)(1-y/\xi)}.
$$
Note that the factors $1-x/\xi,1-y/\xi,1-\xi/x$ and $1-1/y\xi$
correspond to the corners of $D$ which touch the unbounded regions.

\head 3. Potential functions \endhead

Curious to say,
we can always construct a {\it potential function} for the hyperbolicity equations and $\eta$
combinatorially by using Euler's dilogarithm function
$$
\text{Li}_2(z)=-\int_0^z{\log(1-w)\over w}dw,
$$
where we remark that the volume of a tetrahedron with modulus $z$ is given by
$$
D(z)=\text{Im}\,\text{Li}_2(z)+\log|z|\arg(1-z).
$$

\subsubhead 3.1 Neumann-Zagier's functions \endsubsubhead
In fact, we define $V(x,y,\xi)$ by
$$
-\text{Li}_2(1/y\xi)+\text{Li}_2(y/\xi)-\text{Li}_2(y/x)
+\text{Li}_2(\xi/x)+\text{Li}_2(x/\xi)+\log\xi\log{x^2\over y^2\xi^6}-{\pi^2\over6},
$$
the principal part of which is nothing but the sum of dilogarithm functions
associated to the corners of the graph as shown in Figure 7.

\figure
\epsfxsize 2in\epsfbox{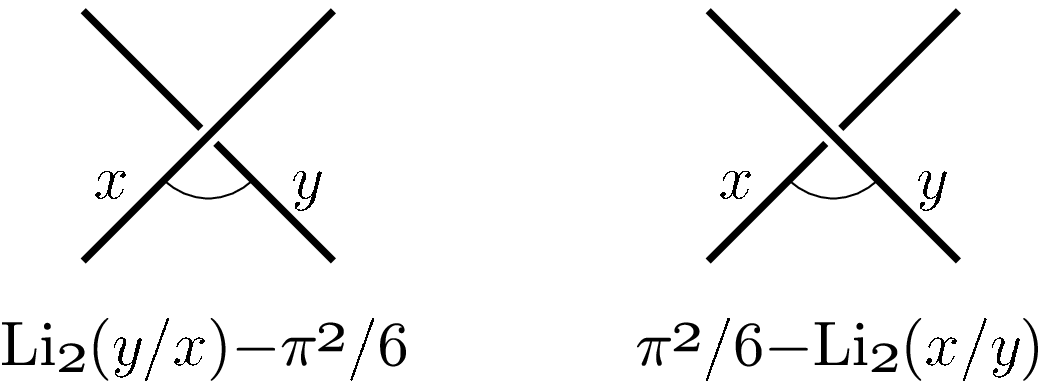}
\endfigure

Then, we have
$$
x{\partial V\over\partial x}=\log{\xi^2(1-\xi/x)\over(1-y/x)(1-x/\xi)},
\quad
y{\partial V\over\partial y}=\log{1-y/x\over\xi^2(1-y/\xi)(1-1/y\xi)},
$$
both of which vanish on $\Cal D$, and
$$
\align
\xi{\partial V\over\partial\xi}&=\log{x^2(1-x/\xi)(1-y/\xi)\over y^2\xi^{12}(1-\xi/x)(1-1/y\xi)}\cr
&=\log\left\{{x\over y\xi^6}\cdot{1\over1-1/y\xi}\right\}^2
-x{\partial V\over\partial x}-y{\partial V\over\partial y}\cr
&=\log\left\{{x\over y\xi^6}\cdot{(1-x/\xi)(1-y/\xi)\over 1-\xi/x}\right\}^2
+x{\partial V\over\partial x}+y{\partial V\over\partial y},
\endalign
$$
that is,
$$
\xi{\partial V\over\partial\xi}=-\log\eta^2
$$
on $\Cal D$,
which shows $V(x,y,e^u)$ coincides with $\Phi(u)$ given in [3,~Theorem 3].

\subsubhead 3.2 Dehn fillings \endsubsubhead
Furthermore, for a {\it slope} $\alpha\in\Bbb Q$, we put
$$
V_\alpha(x,y,\xi)=V(x,y,\xi)+{\log\xi(2\pi\sqrt{-1}-p\log\xi)\over q},
$$
where $p,q\in\Bbb Z$ denote the numerator and the denominator of $\alpha$.
Then, we have
$$
\xi{\partial V_\alpha\over\partial\xi}
=\xi{\partial V\over\partial\xi}+{2\pi\sqrt{-1}-p\log\xi^2\over q}
={2\pi\sqrt{-1}-p\log\xi^2-q\log\eta^2\over q},
$$
and so a solution $(x_\alpha,y_\alpha,\xi_\alpha)$ to the equations
$$
dV_\alpha(x,y,\xi)=0
$$
determines the complete hyperbolic structure of the closed 3-manifold $M_\alpha$
obtained from $M$ by $\alpha$ Dehn filling.
Note that,
by choosing $r,s\in\Bbb Z$ such that $ps-qr=1$,
we can compute the logarithm of the eigenvalue of the holonomy representation
of the core geodesic $\gamma_\alpha$ of $M_\alpha$
which is related to the {\it length} and the {\it torsion} of $\gamma_\alpha$ as follows,
see [3, Lemma 4.2].
$$
\log\xi^r\eta^s={s\pi\sqrt{-1}-\log\xi\over q}
={\text{length}(\gamma_\alpha)+\sqrt{-1}\cdot\text{torsion}(\gamma_\alpha)\over2}.
$$

\subhead Volumes and Chern-Simons invariants \endsubhead

\subsubhead 3.3 Yoshida's functions \endsubsubhead
As in [4], we can observe
$$
\align
\text{Im}\,V_\alpha(x,y,\xi)=&-D(1/y\xi)+D(y/\xi)-D(y/x)+D(\xi/x)+D(x/\xi)\cr
&+\log|x|\cdot\text{Im}\,x{\partial V_\alpha\over\partial x}
+\log|y|\cdot\text{Im}\,y{\partial V_\alpha\over\partial y}
+\log|\xi|\cdot\text{Im}\,\xi{\partial V_\alpha\over\partial\xi},
\endalign
$$
and so
$$
\text{Im}\,V_\alpha(x_\alpha,y_\alpha,\xi_\alpha)=\text{vol}(M_\alpha).
$$
To detect $\text{Re}\,V_\alpha(x_\alpha,y_\alpha,\xi_\alpha)$,
we shall consider
$$
R(x,y,\xi)=-R(1/y\xi)+R(y/\xi)-R(y/x)+R(\xi/x)+R(x/\xi)-{\pi^2\over6},
$$
where $R(z)$ denotes Roger's dilogarithm function defined by
$$
R(z)=\text{Li}_2(z)+\log z\log(1-z)/2.
$$
Then, $R(x,y,\xi)$ can be expressed as
$$
\align
&-\text{Li}_2(1/y\xi)+\text{Li}_2(y/\xi)-\text{Li}_2(y/x)+\text{Li}_2(\xi/x)+\text{Li}_2(x/\xi)\cr
&-{\log x\over2}\!\left(x{\partial V\over\partial x}-\log\xi^2\right)
-{\log y\over2}\!\left(y{\partial V\over\partial y}+\log\xi^2\right)
-{\log\xi\over2}\!\left(\xi{\partial V\over\partial\xi}-\log{x^2\over y^2\xi^{12}}
\right)\!,
\endalign
$$
and so $R(x,y,\xi)$ agrees with
$$
V(x,y;\xi)+\log\xi\,\log\eta
$$
on $\Cal D$ and with
$$
V_\alpha(x,y,\xi)-{\log\xi(2\pi\sqrt{-1}-p\log\xi)\over q}+\log\xi\,\log\eta
=V_\alpha(x,y,\xi)-{\pi\sqrt{-1}\cdot\log\xi\over q}
$$
at $(x_\alpha,y_\alpha,\xi_\alpha)\in\Cal D$.
Therefore, we have
$$
\align
R(x_\alpha,y_\alpha,\xi_\alpha)
&=V_\alpha(x_\alpha,y_\alpha,\xi_\alpha)
-{s\pi^2+\pi\sqrt{-1}\cdot\log\xi_\alpha\over q}+{s\pi^2\over q}\cr
=V_\alpha&(x_\alpha,y_\alpha,\xi_\alpha)
+{\pi\sqrt{-1}\over2}\cdot
\{\text{length}(\gamma_\alpha)
+\sqrt{-1}\cdot\text{torsion}(\gamma_\alpha)\}+{s\pi^2\over q}.
\endalign
$$
In particular,
$$
\align
\text{Im}\,{2\over\pi}\cdot R(x_\alpha,y_\alpha,\xi_\alpha)
&=\text{Im}\,{2\over\pi}\cdot V_\alpha(x_\alpha,y_\alpha,\xi_\alpha)+{2\log|\xi_\alpha|\over q}\cr
&={2\over\pi}\cdot\text{vol}(M_\alpha)+\text{length}(\gamma_\alpha),
\endalign
$$
which shows that, up to a pure imaginary constant,
$$
{2\over\pi\sqrt{-1}}R(x,y,e^u)
$$
must coincide with $2\pi f(u)$ of [3,~Theorem 2], and that
$$
\text{Re}\,{2\over\pi}\cdot R(x_\alpha,y_\alpha,\xi_\alpha)
=\text{Re}\,{2\over\pi}\cdot\left\{V_\alpha(x_\alpha,y_\alpha,\xi_\alpha)
+{s\pi^2\over q}\right\}-\text{torsion}(\gamma_\alpha)
$$
must coincide with $-4\pi CS(M_\alpha)-\text{torsion}(\gamma_\alpha)$.
Consequently, up to some constant which is independent of $\alpha$, we have
$$
\text{Re}\left\{V_\alpha(x_\alpha,y_\alpha,\xi_\alpha)+{s\pi^2\over q}\right\}=-2\pi^2CS(M_\alpha).
$$

\head 4. Concluding remarks \endhead

We redefine $V_\alpha(x,y,\xi)$ as follows.
$$
V_\alpha(x,y,\xi)=V(x,y,\xi)+{\log\xi(2\pi\sqrt{-1}-p\log\xi)+s\pi^2\over q}.
$$
Then, $dV_\alpha(x,y,\xi)=0$ gives the hyperbolicity equations for $M_\alpha$, and
$$
V_\alpha(x_\alpha,y_\alpha,\xi_\alpha)=-2\pi^2CS(M_\alpha)+\text{vol}(M_\alpha)\sqrt{-1}
$$
up to a real constant,
where $(x_\alpha,y_\alpha,\xi_\alpha)$ is a solution to the equations above.

We finally remark that
such a construction always works, even for a {\it link},
and the analytic functions in [3, 5] are now combinatorially constructed
up to a constant.
For the figure-eight knot and $\alpha\in\Bbb Z$,
our potential function coincides with the function in [1]
which appears in the \lq\lq optimistic" limit of the quantum SU(2) invariants of $M_\alpha$.

\references  
\Addresses\recd
\bye